\documentclass[10pt,a4paper]{article}
\usepackage[T1]{fontenc}
\usepackage{amssymb, amsfonts}
\usepackage{amsmath}
\usepackage{physics}
\usepackage{tikz-cd}
\usepackage{enumitem}
\usepackage{comment}

\newcommand{\R}{\mathbb{R}}

\newtheorem{thm}{Theorem}
\newtheorem{lem}{Lemma}
\newtheorem{prop}{Proposition}
\newtheorem{cor}{Corollary}
\newtheorem{defn}{Definition}

\usepackage[backend=bibtex,style=numeric]{biblatex}

\title{A Compactness Condition for the Theorem of Sums in a Class of Non-Uniformly Elliptic PDEs}
\author{Daniel Maienshein}
\date{}

\begin{document}
	\maketitle

	\begin{abstract}

	In the theory of viscosity solutions for second-order, degenerate elliptic PDEs, the Ishii-Lions method is a commonly used strategy, and the theorem of sums is the main analytical tool. As noted by Porretta and Priola, uniformly elliptic PDEs admit a version of the theorem of sums without the squared term due to a compactness argument. Here, we introduce a larger class of PDEs for which the compactness argument holds.

	\end{abstract}

	\section{Notation and Definitions}
	
	We  work in Euclidean space of dimension $N$, and we denote by $\mathbb{S}(N)$ the space of symmetric $N \times N$ matrices with real entries. We  use L\"{o}wner's partial order on $\mathbb{S}(N)$, defined by $A \leq B \iff \langle Ax, x\rangle \leq \langle Bx, x\rangle$ for all $x \in \R^N$. Any symmetric $N\times N$ matrix $X$ possesses $N$ real eigenvalues, and we order them from least to greatest: $\lambda_1(X) \leq \lambda_2(X) \leq \cdots \leq \lambda_N(X)$, where $\lambda_i(X)$ are the eigenvalues of $X$. We also work with the operator norm of a symmetric matrix. That is, we write $\displaystyle \norm{X} = \max_j|\lambda_j(X)|$. We also use the notation $e_i$ for the $i$\textsuperscript{th} standard basis vector in $\R^N$. Moreover, by $\text{diag}(d_1, d_2, \dots, d_N)$ we mean the $N\times N$ diagonal matrix with diagonal entries $d_1, d_2, \dots, d_N$. Finally, suppose that $\mathcal{D} \subseteq \R^N\times \R\times \R^N \times \mathbb{S}(N)$. We call a function $F : \mathcal{D} \to \R$ degenerate elliptic if for every $\omega \in \R^N \times \R\times \R^N$, whenever $(\omega ,X), (\omega ,Y) \in \mathcal{D}$, we have that $X \leq Y \implies F(\omega , X) \geq F(\omega ,  Y)$.

	\section{Introduction}

	Let $\Omega$ be domain in $\R^N$. We work in the generality of fully nonlinear, second-order degenerate elliptic PDEs defined on $\Omega$. Formally, such PDEs can be written as $F(x,u,Du,D^2u) = 0$ for some  degenerate elliptic function $F: \mathcal{D}\to \R$, where $\mathcal{D}$ is a subset of $\Omega\times \R\times \R^N \times \mathbb{S}(N)$. We do not demand  $\mathcal{D} = \Omega\times \R\times \R^N \times \mathbb{S}(N)$, nor do we require $F$ to be continuous or proper. However, \textit{we shall always assume that $F$ is continuous with respect to its matrix entry}.
	
	 Now suppose that $u$ is a viscosity subsolution and $v$ is a viscosity supersolution to the PDE defined by such an $F$. If there exists a $C^2$ function $z(x,y)$ such that $u(x)-v(y) - z(x,y)$ attains a local maximum at some point $( \hat{x}, \hat{y})$ interior to $\Omega\times\Omega$, then the theorem of sums (as presented in \cite[Theorem 3.2]{usersguide}) yields, for all  $\epsilon >0$, the existence of $X_\epsilon,Y_\epsilon \in \mathbb{S}(N)$ such that $ (D_xz( \hat{x}, \hat{y}), X_\epsilon) \in \overline{\mathcal{J}^{2,+}}(u,\hat{x}),\; (-D_yz( \hat{x}, \hat{y}), Y_\epsilon) \in \overline{\mathcal{J}^{2,-}}(v,\hat{y}) $, and
	
	\begin{equation}\label{matrixinequality}
		-\left(\frac{1}{\epsilon} + \norm{A}\right)I \leq  	\begin{bmatrix}
			X_\epsilon & 0 \\
			0 & -Y_\epsilon
		\end{bmatrix} \leq A + \epsilon A^2, \quad \text{ where } A := D^2 z( \hat{x}, \hat{y}).
	\end{equation} 
	
	We shall refer to ``there exists $X_\epsilon, Y_\epsilon\dots$'' along with  inequalities (\ref{matrixinequality}) as ``the conclusion of the theorem of sums.''  (For the definition of jets $\mathcal{J}^{2,\pm}$, we refer the reader to any introduction on viscosity solutions; for example, \cite{katzourakis}.) Compare (\ref{matrixinequality}) to the simpler inequality in the case where $u$ and $v$ are classical (smooth) sub- and super-solutions. There, the condition of a local maximum at $(\hat{x}, \hat{y})$ implies that 
	
	\begin{equation}\label{matrixinequality-2}
		\begin{bmatrix}
			D_{xx}u(\hat{x}) & 0 \\
			0 & -D_{yy}v(\hat{y})
		\end{bmatrix} \leq A, \quad \text{ where } A := D^2 z( \hat{x}, \hat{y}). 
	\end{equation}
	
	Under what structural assumptions on the function $F$ can we improve (\ref{matrixinequality}) to something similar to (\ref{matrixinequality-2}), where the squared term $\epsilon A^2$ has been removed? This question was motivated by a remark by Porretta and Priola in \cite{PP2013}, where they note that uniform ellipticity of $F$ (to be made precise below) is sufficient to get a theorem of sums without  $\epsilon A^2$. They seem to have been the first to make this observation explicit, but these results were already implicitly known in \cite{ishii1990viscosity}, for instance.\footnote{The reader may also recall inequalities such as (3.24) in Remark 3.8 of \cite{usersguide}. While it is true that there is no $\epsilon$ dependence in those matrix inequalities, dependence on $A^2$ has not been removed.} Our new contribution is to introduce a larger class of PDEs, containing non-uniformly elliptic ones, for which we can still recover a theorem of sums without the  $\epsilon A^2$ term.

		\section{Class U -- Uniformly Elliptic Equations}

		Suppose that $\Omega$ is a domain in $\R^N$, $\mathcal{D} \subseteq \Omega\times \R\times \R^N \times \mathbb{S}(N)$, and $F: \mathcal{D} \to \R$. There are various ways of defining ``uniform ellipticity'' for $F$, so for simplicity we shall focus only on the following version of it, which can be found with minor modifications\footnote{We have formulated their condition into one without a time variable and one that explicitly includes the domain of $F$.} in \cite{PP2013}.
		
		\begin{defn}
			Let us say that $F$ belongs to \textbf{Class U} if it is continuous with repsect to its matrix entry, and there exists $\lambda > 0$ and a locally bounded function $H$ such that for all $\omega \in \Omega \times \R \times \R^N$, it is the case that 
			$$\displaystyle B \leq M \text{ in } \mathbb{S}(N) \implies  F(\omega,B) - F(\omega,M) \geq  \lambda \tr(M - B) + H(\omega),$$ provided that  $(\omega,B), (\omega,M) \in \mathcal{D}$. 
		\end{defn}

		In \cite{PP2013}, the authors sketch an argument that whenever $F$ is of Class U and the conclusion of the theorem of sums holds, then the $\epsilon A^2$ term can be removed by a compactness argument. Functions in Class U include $F(X) = -\tr(X)$, leading to the Laplace equation.
		
		We remark that uniform ellipticity is commonly defined using two inequalities rather than one in the literature. See Example 1.12 in \cite{usersguide} or Definition 2.1 in \cite{caffarellicabre1995}, the latter of which is considered only in the case $F = F(x, D^2u)$. However, Condition (4.1) of \cite{jensen1989uniqueness} defines a type of uniform ellipticity with only one inequality, similar to our Class U. The upshot is that no matter which definition of uniform ellipticity is chosen, the removal of $\epsilon A^2$ is possible.

	\section{Class M and The Main Result}
	
	While uniform ellipticity guarantees the removal of $\epsilon A^2$ in the theorem of sums, it is a stronger condition than necessary. In this section, we show that a much larger class of PDEs admits a theorem of sums without $\epsilon A^2$.  As usual, we begin with the assumption that  $F = F(x,u,Du,D^2u)$ is a given degenerate elliptic function defined on  $\mathcal{D} \subseteq \Omega \times \R \times \R^N\times \mathbb{S}(N)$, where $\Omega \subseteq \R^N$. 
	
	\begin{defn}\label{defn-ClassM}
		We shall say that $F$ belongs to \textbf{Class M} if it is continuous in its matrix entry, and given $\omega \in \Omega\times\R\times\R^N$ and $\mathcal{S}_i \subseteq \mathbb{S}(N)$, $i = 1,2$, such that $\{\omega\}\times\mathcal{S}_i \subseteq \mathcal{D}$, then there exist $g_i : \R \times \mathcal{S}_i\to \R$, $i = 1, 2$, (possibly depending on $\omega$) with the following properties:
		
		\begin{enumerate}
			\item For each fixed $M_0 \in \mathcal{S}$, the mappings $t \mapsto g_i(t, M_0)$ are increasing bijections from $\R$ to $\R$. Let us denote the inverses by $ s \mapsto g_i(-, M_0)^{-1}(s)$.
			
			\item The mappings $M \mapsto g_i(-, M)^{-1}(0)$ are continuous on $\mathcal{S}_i$.
	
			\item  $X \leq M \implies -F(\omega,X) \leq g_1(\lambda_1(X), M)$ whenever $(\omega,X) \in \mathcal{D}$, and 
			\item  $-Y \leq M \implies -F(\omega,Y) \geq g_2(\lambda_N(Y), M)$ whenever $(\omega,Y) \in \mathcal{D}$.
		\end{enumerate}
		
	\end{defn}

	We hope to convince the reader that Class M is actually very natural. To begin with, every uniformly elliptic function is of Class M.

	\begin{lem}\label{ClassU-ClassM} Class U is a proper subset of Class M. In particular, if $\lambda > 0$ and $H$ verify membership of $F : \mathcal{D} \to \R$ in Class U, then for $(\omega, M) \in \mathcal{D}$, the following two functions verify membership of $F$ in Class M:  $g_1(t,M) = \lambda t - \lambda \cdot \lambda_1(M) - H(\omega) - F(\omega,M)$, and $g_2(t,M) = \lambda t - \lambda\cdot \lambda_1(M) + H(\omega) - F(\omega,-M)$.  
	\end{lem}
	
	\textbf{Proof.} The fact that there are functions in Class M which are not in Class U is an immediate corollary of Proposition \ref{p-Laplace-example} below, where the $p$-Laplacian is an example for $p \in (1,2)$. Suppose now that $F$ belongs to Class U with $\lambda > 0$ and locally bounded function $H$. Fix $\omega \in \Omega\times\R\times\R^N$ and $M \in \mathcal{S}_1 = \mathcal{S}_2 = \mathbb{S}(N)$, and define $g_1,g_2$ as in the Lemma statement. First, we note that since $\lambda >0$, each $g_i$ is an increasing bijection in the $t$ variable. Moreover, $g_i(-, M)^{-1}(0)$ are continuous in $M$. For instance: $$g_1(-, M)^{-1}(0) = \left(\frac{H(\omega) + F(\omega,M) }{\lambda}+ \lambda_1(M)\right).$$
	
  Therefore, continuity with respect to $M$ follows by continuity of $F$ in the matrix entry along with continuity of eigenvalues.
	
	 If $X \leq M$, then $F(\omega,X) - F(\omega,M) \geq  \lambda \tr(M - X) + H(\omega).$ Hence:
	\begin{align*}
		-F(\omega,X) & \leq  \lambda \tr(X) - \lambda \tr (M) - H(\omega) - F(\omega,M)\\
		&= \lambda \cdot\sum_{j=1}^N \left(\lambda_j(X) -\lambda_j(M) \right) - H(\omega) -F(\omega,M)\\
		&\leq \lambda \cdot (\lambda_1(X) - \lambda_1(M)) - H(\omega) -F(\omega,M)\\
		&= g_1(\lambda_1(X), M) 
	\end{align*}
	
	Similarly, if $-Y \leq M$, then one can show that
	$-F(\omega, Y) \geq g_2(\lambda_N(Y), M)$. 
	\vspace{0.2in}

	We now come to the main result.
	
	\begin{thm}\label{mainthm}
		Suppose that $u$ is a viscosity subsolution and $v$ is a viscosity supersolution to a degenerate elliptic PDE $F(x,u,Du,D^2u) = 0$, where $F : \mathcal{D} \to \R$ belongs to Class M. Suppose that the conclusion of the theorem of sums holds, with $p := D_xz(\hat{x},\hat{y})$, $q := -D_yz(\hat{x},\hat{y})$, and $A := D^2z(\hat{x},\hat{y}).$ If for all $\epsilon >0$ we have $(\hat{x}, u(\hat{x}),p,X_\epsilon), (\hat{y}, v(\hat{y}), q, Y_\epsilon) \in \mathcal{D}$, then there exist $X,Y \in \mathbb{S}(N)$ such that $(p, X) \in \overline{{\mathcal{J}}^{2,+}}(u,\hat{x}), (q, Y) \in \overline{{\mathcal{J}}^{2,-}}(v,\hat{y}) $, and
		\begin{equation*}
			\begin{bmatrix}
				X & 0 \\
				0 & -Y
			\end{bmatrix} \leq A,
		\end{equation*} 	
	
	\vspace{0.1in}

	Moreover, using the notation of Definition \ref{defn-ClassM}, there exist $g_1,g_2$ yielding the lower bounds  $$[\,g_1(-, D_{xx}z(\hat{x},\hat{y}))^{-1}(0)\,]I   \leq X,$$ and $$-[\,g_2(-,D_{yy}z(\hat{x},\hat{y}))]^{-1}(0)\,]I \leq -Y.$$
	
	Finally, $g_1$ may depend only on $F, \hat{x}, u,$ and $p$, while $g_2$ may depend only on $F,\hat{y}, v,$ and $q$.

\end{thm}	

\vspace{0.1in}

Before proving Theorem \ref{mainthm}, let us make two remarks. First, the assumptions $(\hat{x}, u(\hat{x}),p,X_\epsilon), (\hat{y}, v(\hat{y}), q, Y_\epsilon) \in \mathcal{D}$ are necessary for Theorem \ref{mainthm} to be applicable when $F$ is not continuous or defined everywhere. For example, equations of the $p$-Laplace type may have singularities when $Du = 0$. We shall elaborate on this point in the Examples section below.

Second, the lower bounds depend on the functions $g_1$ and $g_2$, rendering them difficult to interpret as is. Hence, it is worthwhile to see what they become in the more concrete setting of uniformly elliptic equations. The following is a corollary of Theorem \ref{mainthm} along with Lemmas \ref{ClassU-ClassM} and \ref{upperbound}.

\begin{cor}\label{maintheorem-ClassU}
	Given the setup of Theorem \ref{mainthm}, if $F$ is of Class U with $\lambda > 0 $ and locally bounded function $H$, then there exist $X,Y \in \mathbb{S}(N)$ such that $(p, X) \in \overline{{\mathcal{J}}^{2,+}}(u,\hat{x}), (q, Y) \in \overline{{\mathcal{J}}^{2,-}}(v,\hat{y})$, and
	\begin{equation*}
		\begin{bmatrix}
			X & 0 \\
			0 & -Y
		\end{bmatrix} \leq A,
	\end{equation*}

	where $p = D_xz(\hat{x},\hat{y})$, $q = -D_yz(\hat{x},\hat{y})$, and $A = D^2z(\hat{x},\hat{y}).$ \vspace{0.1in}
	
	Moreover, we have the lower bounds  $$\left(\frac{H(\hat{x}, u(\hat{x}), p) + F(\hat{x}, u(\hat{x}), p,D_{xx}z(\hat{x},\hat{y}) ) }{\lambda}+ \lambda_1(D_{xx}z(\hat{x},\hat{y}))\right) I   \leq X,$$ and $$ \left(\frac{H(\hat{y}, v(\hat{y}), q) - F(\hat{y}, v(\hat{y}), q,-D_{yy}z(\hat{x},\hat{y}) ) }{\lambda}- \lambda_N(-D_{yy}z(\hat{x},\hat{y}))\right) I \leq -Y.$$

\end{cor}

We now progress toward the proof of Theorem \ref{mainthm}. The first lemma we need is mentioned in \cite{PP2013}. We present a slightly more flexible version here. It provides us with a uniform upper bound on matrices $X_\epsilon$ and $-Y_\epsilon$ appearing in the conclusion of the theorem sums. 

\begin{lem}\label{upperbound}
	Suppose that $A = \begin{bmatrix}
		E & B \\
		B^\top & D
	\end{bmatrix}$ and fix $\epsilon_0 > 0$. If for all $\epsilon \in (0,\epsilon_0)$ there exist symmetric matrices $X_\epsilon, -Y_\epsilon$ satisfying	
	$$\begin{bmatrix}
		X_\epsilon & 0 \\
		0 & -Y_\epsilon
	\end{bmatrix} \leq A + \epsilon A^2,\quad$$  then $X_\epsilon \leq E + \epsilon_0(E^2 + B B^\top)$ and $-Y_\epsilon \leq D + \epsilon_0(D^2 + B^\top B)$ for all $\epsilon \in (0,\epsilon_0)$. 
	
\end{lem}

\textbf{Proof.} Note that $A^2$ is symmetric and positive semi-definite, and $A^2 = \begin{bmatrix}
	E^2 + BB^\top & * \\
	* & D^2 + B^\top B
\end{bmatrix}.$  The result follows from computing quadratic forms with vectors $\begin{bmatrix}
	\xi \\
	0
\end{bmatrix} \text{ and } \begin{bmatrix}
	0\\
	\xi
\end{bmatrix}$ in $\R^{2N}$ and using that $\epsilon < \epsilon_0$.\newline

Note that Lemma \ref{upperbound} applies to any degenerate elliptic $F$. However, we shall also need a uniform lower bound on $X_\epsilon$ and $-Y_\epsilon$. First, we should remark that the theorem of sums actually does provide us with a lower bound: $$-\left(\frac{1}{\epsilon} + \norm{A}\right)I \leq  	\begin{bmatrix}
	X_\epsilon & 0 \\
	0 & -Y_\epsilon
\end{bmatrix}.$$ Unfortunately, this lower bound does not provide the compactness we need, because it is not uniform in $\epsilon$. Therefore, most of the work behind Theorem \ref{mainthm} goes into establishing lower bounds which do not depend on $\epsilon$. \vspace{0.1 in}

The other lemma we need is stability of jet closures, a viscosity theory fact which is known but not explicitly written down in the literature as far as the author is aware.\footnote{Similar statements can be found in  \cite{usersguide} or \cite{katzourakis} under the name of closure of ``slice sets.'' However, these only apply--strictly speaking--to genuine jets and not jet \textit{closures}.} For sake of completeness, we include a proof here. 

\begin{lem}\label{jetclosures-are-closed}
	Suppose that there exists $x \in \R^N$, $p \in \R^N$ and a sequence of symmetric $N \times N$ matrices $X_m$ such that for all $m$, we have $(p,X_m) \in \overline{\mathcal{J}^{2,+}}(u,x)$. If $X_m \to X$ as $m \to \infty$, then $(p,X) \in \overline{\mathcal{J}^{2,+}}(u,x)$.
\end{lem}

\textbf{Proof.} The hypothesis that $(p,X_m) \in \overline{\mathcal{J}^{2,+}}(u,x)$ implies, by definition, that for each $m$, there exists $(x^k, u(x^k), p^k, X_m^k)$ converging to $(x,u(x),p,X_m)$  and $(p^k, X_m^k) \in \mathcal{J}^{2,+}(u,x^k)$ for all $k$. We thus obtain a doubly-indexed sequence of matrices $X_m^k$ with $X_m^k \to X_m$ for each $m$, and $X_m \to X$. There is a subsequence $X_m^{k(m)}$ converging to $X$ as $m \to \infty$. To see this, note that for each entry $\xi_m^k$ of $X_m^k$,  $\xi_m$ of $X_m$, and $\xi$ of $X$, we have $|\xi_m^{k(m)} - \xi| \leq |\xi_m^{k(m)} - \xi_m|+|\xi_m -\xi|$ for $k(m)$ chosen such that $|\xi_m^{k(m)} - \xi_m| < m^{-1}$. Therefore, we have $(x^{k(m)}, u(x^{k(m)}), p^{k(m)}, X_m^{k(m)})$ converging to $(x,u(x),p,X)$  and $(p^{k(m)}, X_m^{k(m)}) \in \mathcal{J}^{2,+}(u,x^{k(m)})$ for all $m$. Thus by definition, $(p,X)  \in \overline{\mathcal{J}^{2,+}}(u,x)$.

\vspace{0.2 in}

\textbf{Proof of Theorem \ref{mainthm}.} Temporarily fix $\epsilon_0 > 0$, and suppose that $0 < \epsilon < \epsilon_0.$ Let $E = D_{xx}z(\hat{x},\hat{y})$, $\tilde{E} = (D_{xx}z(\hat{x},\hat{y}))^2 + D_{xy}z(\hat{x},\hat{y})D_{yx}z(\hat{x},\hat{y})$. By Lemma \ref{upperbound} and since $\tilde{E} \geq 0$, we have $X_\epsilon \leq E + \epsilon\tilde{E}\leq E + \epsilon_0\tilde{E}$. Write $\mathcal{S}_1 = \{ E + \epsilon\tilde{E} : \epsilon\geq  0\}$. Since $F$ is of Class $M$, there exists $g_1 : \R\times\mathcal{S}_1 \to \R$, possibly depending on $F, \hat{x}, u, p$ satisfying conditions 1 through 4 of Definition \ref{defn-ClassM}, such that

$$ -F(\hat{x}, u(\hat{x}), p, X_\epsilon) \leq  g_{1}(\lambda_{1}(X_\epsilon), E + \epsilon_0\tilde{E}).$$

Since $u$ is a viscosity subsolution and $(\hat{x},u(\hat{x}),p,X_\epsilon) \in \mathcal{D},$ we know that $F(\hat{x},u(\hat{x}),p,X_\epsilon) \leq 0.$ By condition 1 of Definition \ref{defn-ClassM}, we deduce that
$$g_{1}(-, E + \epsilon_0\tilde{E})^{-1}(0)\leq  \lambda_{1}(X_\epsilon) \leq \lambda_{N}(X_\epsilon) \leq \lambda_N(E + \epsilon_0\tilde{E}).$$ Therefore, matrices $X_\epsilon$ are bounded in the operator norm, uniformly in $\epsilon \in (0,\epsilon_0)$. Thus there exists a subsequence in $\epsilon$, along which $X_\epsilon$ converge to some symmetric matrix $X$. Taking the limit as $\epsilon \to 0$ and using continuity of eigenvalues, we obtain $g_{1}(-, E + \epsilon_0\tilde{E})^{-1}(0) \leq  \lambda_{1}(X).$

We now repeat the same procedure to bound $Y_\epsilon$. Suppose again that $0 < \epsilon < \epsilon_0.$ Let $D = D_{yy}z(\hat{x},\hat{y})$, $\tilde{D} = (D_{yy}z(\hat{x},\hat{y}))^2 + D_{yx}z(\hat{x},\hat{y})D_{xy}z(\hat{x},\hat{y})$. By Lemma \ref{upperbound} and since $\tilde{D} \geq 0$, we have $-Y_\epsilon \leq D + \epsilon\tilde{D}\leq D + \epsilon_0\tilde{D}$. Write $\mathcal{S}_2 = \{  D + \epsilon\tilde{D} : \epsilon \geq 0\}$. Since $F$ is of Class $M$, there exists $g_2 : \R\times\mathcal{S}_2 \to \R$, possibly depending on $F, \hat{y}, v, q$ satisfying conditions 1 through 4 of Definition \ref{defn-ClassM}, such that 

$$-F(\hat{y}, v(\hat{y}), q, Y_\epsilon)  \geq  g_{2}(\lambda_{N}(Y_\epsilon),D +\epsilon_0\tilde{D} ).$$

Since $v$ is a viscosity subsolution and $(\hat{y},v(\hat{y}),q,Y_\epsilon) \in \mathcal{D},$  we know that $F(\hat{y},v(\hat{y}),q,Y_\epsilon) \geq 0.$ By condition 1 of Definition \ref{defn-ClassM}, we deduce that
$$g_{2}(-,D + \epsilon_0\tilde{D} )^{-1}(0) \geq  \lambda_{N}(Y_\epsilon) \geq \lambda_{1}(Y_\epsilon) \geq \lambda_1(-D - \epsilon_0\tilde{D}).$$ Therefore, matrices $Y_\epsilon$ are bounded in the operator norm, uniformly in $\epsilon \in (0,\epsilon_0)$. Thus there exists a subsequence in $\epsilon$, along which $Y_\epsilon$ converge to some symmetric matrix $Y$; without loss of generality we can assume that $X \to X_\epsilon$ along this same subsequence. Taking the limit as $\epsilon \to 0$, we obtain $g_{2}(-, D + \epsilon_0\tilde{D})^{-1}(0) \geq \lambda_{N}(Y).$  

Noting that we have $\displaystyle \begin{bmatrix}
	X& 0 \\
	0 & -Y
\end{bmatrix} \leq A + \epsilon_0A^2$, we take the limit as $\epsilon_0 \to 0$ to obtain the  upper bound $\displaystyle \begin{bmatrix}
	X & 0 \\
	0 & -Y
\end{bmatrix} \leq A$, as well as the  lower bounds (using the continuity condition 2 of Definition \ref{defn-ClassM} for $g_1$ and $g_2$) $$g_{1}(-, E)^{-1}(0)\leq \lambda_1(X),$$ and $$-g_{2}(-, D)^{-1}(0) \leq -\lambda_N(Y) = \lambda_1(-Y).$$ 

Finally, invoking Lemma \ref{jetclosures-are-closed}, we obtain that $(p, X) \in \overline{{\mathcal{J}}^{2,+}}(u,\hat{x})$ and $(q, Y) \in \overline{{\mathcal{J}}^{2,-}}(v,\hat{y})$. This completes the proof.

\subsection{Examples and Non-Examples}\label{examples-new} In this section we provide,  with proof, examples of functions belonging to (or not belonging to) Class M. Some of these examples already belong to Class U, so the extra information in those cases is to make explicit the lower bounds from Theorem \ref{mainthm} (when the Theorem applies).

\subsubsection{Uniformly Elliptic Linear Equations}
	
 Let $F(x,u, Du, D^2u) = -\tr (a(x) D^2u) + b(x)\cdot Du + c(x) u$, where $a$ is uniformly elliptic: there exists $\theta > 0$ and positive semi-definite matrix-valued function $\sigma(x)$ such that for all $x \in \R^N$, we have $a(x) = \sigma(x) + \theta I$. 
 
 These equations all belong to Class U. One can take $\lambda = \theta$, and $H \equiv 0$. Indeed, whenever $B \leq M$, we have 
\begin{align*}
	F(x, z, \nu, B) - F(x, z, \nu, M) &= \tr(a(x)M) - \tr(\sigma(x)B) \\
	&= \theta\tr(M-B) + \tr(\sigma(x)(M-B)) \\
	&\geq \theta\tr(M-B)
\end{align*} 

We have used the general fact that whenever $A \leq A'$ and $C \geq 0$ in $\mathbb{S}(N)$, then $\tr(C A) \leq \tr(CA')$. Therefore, the lower bounds in Corollary \ref{maintheorem-ClassU} are:

$$\left(\frac{-\tr (a(x) D_{xx}z(\hat{x},\hat{y})) + b(\hat{x})\cdot p + c(\hat{x}) u(\hat{x})}{\theta}+ \lambda_1(D_{xx}z(\hat{x},\hat{y}))\right) I   \leq X,$$ and $$ \left(\frac{ -\tr (a(x) D_{yy}z(\hat{x},\hat{y})) - b(\hat{y})\cdot q - c(\hat{y}) v(\hat{y})}{\theta}- \lambda_N(-D_{yy}z(\hat{x},\hat{y}))\right) I \leq -Y.$$

\subsubsection{The $p$-Laplace Equation}

We consider the $p$Laplace equation $-\text{div}(|\nabla u|^{p-2}\nabla u) = 0$ for $1 \leq p < \infty$. The equation can formally be re-written as $F_p(\nabla u, D^2u) = 0$, where\\ $F_p : (\R^N \setminus\{0\})\times\mathbb{S}(N) \to \R$ is given by 

$$F_p(\nu, X) = -|\nu|^{p-2}\left[\tr X + (p-2)\left\langle X \frac{\nu}{|\nu|}, \frac{\nu}{|\nu|}\right\rangle\right].$$

While $F_p$ is degenerate elliptic for all values of $p$, it is known that $F_p$ does not belong to Class U whenever $p \neq 2$. We briefly verify this fact. For any $\lambda > 0$ and locally bounded $H$, we can always choose $c = c(\lambda) > 0$ so that $\nu := ce_1$ satisfies $|\nu|^{p-2} - \lambda < 0$, and hence we can find $l > 0$ large enough so that, locally, $(|\nu|^{p-2} - \lambda)l < \inf H$. It follows that for $Y = \text{diag}(0, l)$ and $X = \text{diag}(0,0)$, $X\leq Y$. Yet due to the fact that $\nu$ is in the nullspace of $Y-X$, we have $$F_p(\nu, X) - F_p(\nu, Y) = |\nu|^{p-2}\, l < \lambda \tr(Y-X) + H(\nu).$$

On the other hand, we have the following: 

\begin{prop}\label{p-Laplace-example}
	Function $F_p$ belongs to Class M if and only if $p \in (1,\infty)$.
\end{prop}

\textbf{Proof.} Continuity of $F$ with respect to the matrix entry is evident. Let $\nu \in \R^N \setminus\{0\}$, and define $\mathcal{S}_1 = \mathcal{S}_2 = \mathbb{S}(N)$. First consider the case $p \geq 2$. Let $g_1(t, M) = |\nu|^{p-2}\left[t + (N+p-3)\lambda_N(M)\right]$ and $g_2(t,M) = |\nu|^{p-2}\left[t - (N+p-3)\lambda_N(M)\right]$.

Properties 1 and 2 of Definition \ref{defn-ClassM} can easily be checked, and we find that $g_1(-, M)^{-1}(0) = - (N+p-3)\lambda_N(M)$ and $g_2(-, M)^{-1}(0) = (N+p-3)\lambda_N(M)$.

Moreover, whenever $X \leq M$, we have 
\begin{align*}
	-F_p(\nu, X) &\leq |\nu|^{p-2} \left[\lambda_1(X) + \sum_{j=2}^N \lambda_j(X) + (p-2)\lambda_N(X)\right]\\
	&\leq  |\nu|^{p-2} \left[\lambda_1(X) + (N-1)\lambda_N(M) + (p-2)\lambda_N(M)\right]\\
	&=  |\nu|^{p-2} \left[\lambda_1(X) + (N+p-3)\lambda_N(M)\right] \\
	&=g_1(\lambda_1(X),M).
\end{align*}

Similarly, if $-Y \leq M$, then $-F_p(\nu, Y) \geq g_2(\lambda_N(B), M).$

\vspace{0.1in}

Next we tackle the case $p \in (1,2)$. Let $g_1(t, M) = |\nu|^{p-2}\left[(p-1)t + (N-1)\lambda_N(M)\right]$ and $g_2(t,M) = |\nu|^{p-2}\left[(p-1)t - (N-1)\lambda_N(M)\right]$.

Properties 1 and 2 of Definition \ref{defn-ClassM} can be verified, and we find that $\displaystyle g_1(-, M)^{-1}(0) =  -\frac{N-1}{p-1}\lambda_N(M)$ and $\displaystyle g_2(-, M)^{-1}(0) =\frac{N-1}{p-1}\lambda_N(M)$.

Moreover, whenever $X \leq M$, we have 

\begin{align*}
	-F_p(\nu, X) &\leq |\nu|^{p-2} \left[(p-1) \tr(X) - (p-2)\tr(X)  + (p-2)\lambda_1(X)\right]\\
	&\leq |\nu|^{p-2} \left[(p-1)\lambda_1(X) + \sum_{j\geq 2} \lambda_j(X)\right]\\
	&\leq  |\nu|^{p-2}\left[(p-1)\lambda_1(X) + (N-1) \lambda_N(M)\right]\\
	&=g_1(\lambda_1(X),M).
\end{align*}

Similarly, if $-Y \leq M$, then $-F_p(\nu, Y) \geq g_2(\lambda_N(B), M).$ 

Finally, if $p = 1$, then $F_1$ is not of Class M. For a counterexample, we can take $\nu = e_N$, $\mathcal{S}_1 = \mathcal{S}_2 = \mathbb{S}(N)$, and $M = I$. Let $g_i: \R \times \mathcal{S}_i \to \R$ be arbitrary functions satisfying conditions 1 through 4 of Definition \ref{defn-ClassM}. Since $t \mapsto g_1(t, I)$ is unbounded below, there exists $c \leq 0$ such that $g_1(c,I) < N-1$. Define $X = \text{diag}(1, 1, \dots, 1, c)$. Then $X \leq I$, $\lambda_1(X) = c$, and yet we have:	
\begin{align*}
	-F_1(\nu,X) &= |1|^{-1}\left[(N-1 + c) - c \right] > g(c,I),
\end{align*} 

which contradicts condition 3. The proof is complete. \vspace{0.1in}

As a corollary of Proposition \ref{p-Laplace-example}, for $p \in (1, \infty)$ we get the following lower bounds from Theorem \ref{mainthm}:

If $p \geq 2$, then
 $$- (N+p-3)\lambda_N(D_{xx}z(\hat{x},\hat{y}))\,I   \leq X,$$ and $$-(N+p-3)\lambda_N(D_{yy}z(\hat{x},\hat{y}))\,I \leq -Y.$$
 
If $1 < p < 2$, then
 $$-\frac{N-1}{p-1}\lambda_N(D_{xx}z(\hat{x},\hat{y}))\,I   \leq X,$$ and $$-\frac{N-1}{p-1}\lambda_N(D_{yy}z(\hat{x},\hat{y}))\,I \leq -Y.$$

It is therefore not surprising that $F_1$ fails to be of Class M. Indeed, as $p \to 1^+$, the lower bounds blow up. We also note that when $p \to +\infty$, the lower bounds blow up as well, which indicates that the $\infty$-Laplacian should also fail to be in Class M. We will prove this below, but first let us make a remark about the notion of viscosity sub-/super-solution to $p$-Laplace equations for $p \in (1,\infty)$.

 The definition of supersolution (for instance) requires a bit of extra care, because the function $F_p$ is not defined at $\nu = 0$ and indeed has no continuous extension there (for $N \geq 2$). Following \cite{JLMequivalence}, a function $u : \Omega \to (-\infty, \infty]$ is a viscosity supersolution to $F_p(\nabla u, D^2u) = 0$ if and only if  $u$ is lower semicontinuous, not identically $+\infty$, and whenever $\phi \in C^2(\Omega)$ touches $u$ from below at $x_0$, if $\nabla\phi(x_0) \neq 0$, then $F_p(\nabla \phi, D^2\phi) \geq 0$. Note that this is almost identical to a standard definition of viscosity supersolution, except for the fact that we make no assertions when $\nabla \phi$ vanishes. The point is that this definition is the correct one, because it makes viscosity solutions equivalent to $p$-harmonic functions for $p \in (1,\infty)$; see \cite[Corollary 2.6]{JLMequivalence}. 
 
 What is the relevance of these comments to this paper? The domain $\mathcal{D}$ for $F_p$ is $\R^N\setminus\{0\}\times \mathbb{S}(N)$. For Theorem \ref{mainthm} to apply to $F_p$, we need to assume that  $(D_xz(\hat{x},\hat{y}),X_\epsilon), (-D_yz(\hat{x},\hat{y}), Y_\epsilon) \in \mathcal{D}$, so we need $-D_yz(\hat{x},\hat{y})\neq 0$, for example. Of course, all we know in general is that $(-D_yz( \hat{x}, \hat{y}), Y_\epsilon) \in \overline{\mathcal{J}^{2,-}}(u,\hat{y}) $, so it is generally possible that $ D_yz(\hat{x},\hat{y}) =  0$. Hence, the definition of viscosity supersolution tells us nothing about $F_p(-D_yz(\hat{x},\hat{y}), Y_\epsilon)$; indeed $F_p$ is not even defined at $\nu = 0$.  Fortunately, with an appropriate choice of test function $z$, we can force $-D_yz(\hat{x},\hat{y})\neq 0$ and $D_xz(\hat{x},\hat{y}) \neq 0$; this is demonstrated in the proof of \cite[Proposition  3.3]{JLMequivalence}.  

\subsubsection{The homogeneous $p$-Laplace Equation}

 We could also consider the homogeneous $p$-Laplace functions $F_p^H : (\R^N \setminus\{0\})\times\mathbb{S}(N) \to \R$ given by 
 
 $$F_p^H(\nu, X) = -\tr X - (p-2)\left\langle X \frac{\nu}{|\nu|}, \frac{\nu}{|\nu|}\right\rangle.$$
 
 While the functions $F_p^H$ still have no continuous extension at $\nu = 0$ for $N \geq 2$, they can be shown to belong to Class U for $p \in (1,\infty)$. Also, the lower bounds relevant for Theorem \ref{mainthm} are identical with the lower bounds for the (non-homogeneous) $p$-Laplacians, because the factor of $|\nu|^{p-2}$ plays no role when the inverse functions $g_i(-,M)^{-1}$ are evaluated at $0$. Therefore, we shall not discuss these equations any further.

\subsubsection{The $\infty$-Laplace Equation}

The $\infty$-Laplace equation is: $-\Delta_\infty u = 0$, where $\Delta_\infty u = \langle D^2u Du, Du\rangle.$ The associated $F_\infty:\R^n\times\mathbb{S}(N) \to \R$ is therefore $F_\infty(\nu,X) = -\langle X \nu, \nu\rangle.$ We claim that $F_\infty$ is not of Class M. Take $N \geq 2$, $\nu = e_1$,  $M = I$, and $X = \text{diag}(1, 0, \dots, 0, c) $, where it is understood that $X$ has $N - 2$ zero diagonal entries. Note that $X \leq M$. If $g_1$ is given as per the definition of Class M, then there exists $c \in \R$ such that $c \leq 0$ and $g(c, I) < 1$. Since $\lambda_1(X) = c$, it follows that 
$$
-F_\infty(\nu,X) = 1  > g(c, I),$$ which is a contradiction.

\vspace{0.1in}

One can similarly show that the homogeneous version $\displaystyle F_\infty^H(\nu,X) = -\frac{\langle X \nu,\nu\rangle}{\langle \nu, \nu\rangle}$ does not belong to Class M. As a corollary, we have verified the known fact that the $\infty$-Laplace equations (homogeneous or not) are not uniformly elliptic.

\subsubsection{Functions of the Eigenvalues of the Hessian}

 \begin{enumerate}

 	\item (the $k$-Hessian) For $k \in \{1,2,\dots, N\}$ and $x \in \R^N$, define the $k$\textsuperscript{th} elementary symmetric polynomial: $S_k(x) = \displaystyle \sum_{i_1 < i_2 < \cdots < i_k} x_{i_1} x_{i_2}\cdots x_{i_k}$. For $X \in \mathbb{S}(N)$, write $\lambda(X) = (\lambda_1(X), \dots, \lambda_N(X))$, and define $F_k(X) = -S_k(\lambda(X)).$ The expression $F_k(D^2u)$ is called the $k$-Hessian of $u$. The $k$-Hessian generalizes the Monge-Amp\`ere equation $-\det(D^2u) = 0$, obtained by setting $k = N$, and the Laplace equation, by setting $k = 1$.
 	
 	To make sense of viscosity solutions, we need to restrict the domain of $F_k$ to a proper subset $\Sigma_k$ of $\mathbb{S}(N)$ to ensure that $F_k$ is degenerate elliptic. Explicitly, $\Sigma_k = \{X \in\mathbb{S}(N) : \lambda(X) \in \overline{\Gamma}_k\}$, where $\displaystyle \Gamma_k = \bigcap_{j=1}^k\{\lambda \in \R^N : S_j(\lambda) > 0\}$, and $\overline{\Gamma}_k$ denotes the closure of $\Gamma_k$. See \cite{wang2009k} for details. 
 	
 	We claim that if $N \geq 2$ and $2 \leq k \leq N$, then $F_k(X) = -S_k(\lambda(X))$ is not of Class M. Take $M = I$, and define $X_n = \text{diag}(-n, -n, 1, \dots, 1)$ for every $n \in \mathbb{N}$. Note that $X_n \leq I$ and $\lambda_1(X_n) = -n$ for every $n$. By counting powers of $n$, we have for each $k \in \{2, 3, \dots, N\}$: 
 	
 	$$-F_k(X_n) = {N-2 \choose k - 2} n^2 - 2{N-2 \choose k - 1}n + {N-2 \choose k}.$$  
 	
 	If $g_1$ is as in the definition of Class M, then we can choose $n$ sufficiently large so that both $-F_k(X_n)  > 0$ and $g(-n, I) < 0$ hold, and we reach the contradiction $-F_k(X_n)  > g(\lambda_1(X_n), I)$.
 	
 	\item (Sums of Functions of Eigenvalues) 
 	
 	Here we explore another class of equations which generalize the Laplace equation. Suppose that $H : \R \to \R$ is strictly increasing. Define \\$\displaystyle F_H(X) = -\sum_{j=1}^N H(\lambda_j(X))$. Then $F_H$ is degenerate elliptic due to monotonicity of $H$.
 	
 	\begin{prop}
 		Function $F_H$ belongs to Class M if and only if $H$ is unbounded above and below. 
 	\end{prop}

 	 \textbf{Proof.} First suppose that $H$ is unbounded above and below. It follows that $H$ is an increasing bijection, and therefore continuous. So for $M \in \mathbb{S}(N)$, the functions $g_1(t,M) = H(t) + \sum_{j=2}^N H(\lambda_j(M))$ and $g_2(t,M) = H(t) + \sum_{j=1}^{N-1} H(-\lambda_j(M))$ verify conditions 1 and 2 of Definition \ref{defn-ClassM}. Further, whenever $X \leq M$, we have:
 	\begin{align*}
 		-F_H(X) &= \sum_{j=1}^N H(\lambda_j(X))  \\
 		&\leq H(\lambda_1(X)) + \sum_{j=2}^N H(\lambda_j(M))\\
 		&= g_1(\lambda_1(X), M) 
 	\end{align*}
 	
 	Likewise, if $-Y \leq M$, then $-F_H(Y) \geq  g_2(\lambda_N(Y), M) $

 	Next, we suppose that $H$ is bounded above or below; without loss of generality suppose that $H$ is bounded below by some constant $K \in \R$. Let $N = 1$ and $\mathcal{S}_1 = \{-n : n \in \mathbb{N}\}$. Suppose for sake of contradiction that $g_1$ existed as per the definition of Class M. Let $X_n = -n$ for all $n \in \mathbb{N}$, and let $M = 0.$ Since $\lambda_1(X_n) = -n$, membership in Class M implies the following sequence of inequalities:
 	$$K \leq H(-n) = -F_H(X_n) \leq g_1(-n, 0).$$ Taking the limit as $n \to \infty$, and since $g_1$ is an increasing bijection in its first entry, we find that 
 	$$K \leq -\infty,$$ which is a contradiction. 
 	
 	\vspace{0.1in}
 	
 	In case $H$ is unbounded above and below, we get the lower bounds 
 	
 	 $$H^{-1}\left[-\sum_{j=2}^N H(\lambda_j(D_{xx}z(\hat{x},\hat{y})))\right]\,I   \leq X,$$ and $$-H^{-1}\left[-\sum_{j=1}^{N-1} H(-\lambda_j(D_{yy}z(\hat{x},\hat{y})))\right]\,I \leq -Y.$$
 	
 	An interesting example of an $H$ which is unbounded above and below is $H(t) = t^{1/d}$, where $d \geq 3$ is an odd integer and $N \geq 2$. The function $F(X) = -\sum_{j=1}^N [\lambda_j(X)]^{1/d}$ therefore belongs to Class M but not to Class U. Let us prove this. Otherwise, there would exist $\lambda >0$ and a constant $K$ such that  
 	$\displaystyle B \leq M \text{ in } \mathbb{S}(N) \implies  F(B) - F(M) \geq  \lambda \tr(M-B) + K.$ But for $M = 0$ and $X_n = -\text{diag}(n, 1/n, \dots, 1/n)$, we have $X_n \leq M$, and the following sequence of inequalities holds:
 	
 	$$n^{1/d} + \frac{(N-1)}{n^{1/d}} \geq  -\lambda \left(-n - \frac{(N-1)}{n}\right) + K.$$
 	
 	These inequalities imply that $\lambda n \leq n^{1/d} - K + O(n^{-1/d})$ as $n \to \infty$, a contradiction. 
 	
 	\vspace{0.1 in}
 	
 	An example of a bounded strictly increasing $H$ is $H(t) = \arctan(t),$ leading to the Special Lagrangian equation, which is studied in \cite{yuan2020special}.

	 \end{enumerate}

\section{Discussion about Comparison Principles}

Thanks to Theorem \ref{mainthm}, the conclusion of the theorem of sums can be improved when working in Class M. A natural question to ask is whether viscosity solutions to PDEs defined by functions in Class M enjoy a comparison principle (which implies uniqueness). We mean that a degenerate elliptic PDE $F(x,u,Du,D^2u) = 0$ on $\Omega \subseteq \R^N$ satisfies a comparison principle if for all viscosity subsolutions $u$ and for all viscosity supersolutions $v$, whenever $u \leq v$ on $\partial\Omega$, we have $u \leq v$ in $\Omega$. Here, we show that the comparison principle is a completely separate question from membership in Class M. 

First, there are functions in Class M which do not enjoy a compasion principle. A classic example is the function $F: \R^N \times \mathbb{S}(N) \to \R$ defined by $F(\nu,X) = -\tr X - |\nu|^{1/2}$. Then $F$ belongs to Class U (and hence Class M). However, as noted in \cite{kawohl2000comparison}, one can find two distinct classical solutions to $F(Du,D^2u) = 0$ in $\Omega = \{x \in \R^N : |x| < 1\}$ subject to zero boundary data. So uniqueness (and hence a comparison principle) fails. Typically, one has to assume more about the structure of $F$ to prove a comparison principle. An example of a sufficient condition for uniqueness is Lipschitz continuity in the gradient and no explicit dependence on the $x$ variable; see \cite[Theorem 1]{kawohl2000comparison}. 

Next, there are functions in Class $M$ which enjoy a comparison principle. Equations such as $F(X) = -\tr(X)$ immediately come to mind. For a more interesting example, when $N \geq 2$ and $d \geq 3$ is odd, in the Examples section we mentioned  the function $F(X) =-\sum_{j=1}^N [\lambda_j(X)]^{1/d}$. The author is not aware of this $F$ being considered explicitly in the literature; however, it does fall into various known classes of PDEs which enjoy a comparison principle. For example, $F$ satisfies  part (iv) of \cite[Proposition 2.6]{brustad2023comparison}. Indeed, due to strict monotonicity of $s \mapsto s^{1/d}$, we have for each eigenvalue $\lambda_j$ and for all $t > 0$: $(\lambda_j + t)^{1/d} > \lambda_j^{1/d}$, and $(\lambda - t)^{1/d} < \lambda^{1/d}$. Alternatively, one can deduce a comparison principle by showing that $F$ is ``locally strictly elliptic'' as defined by \cite{kawohl2000comparison} and applying their Theorem 1.  

We can also find functions not in Class M for which a comparison principle holds. Consider the Special Lagrangian function $F(X) = -\sum_{j=1}^N \arctan(\lambda_j(X))$. We argued in the Example section that $F$ is not of Class M since the arctangent is bounded. On the other hand, like the last example, this $F$ falls into the class of PDEs which are considered by \cite{brustad2023comparison} or \cite{kawohl2000comparison}, so that a comparison principle holds.

Finally, the $1$-Laplacian is a good example of a function which is neither in Class M nor satisfies a comparison principle. We already showed above that it is not of Class M. Moreover, by \cite[Example 3.6]{sternberg1994generalized}, uniqueness (and hence a comparison principle) does not hold.

\section{Acknowledgments}

The author expresses gratitude to Juan Manfredi, Fausto Ferrari, and Nicol\`o Forcillo for helpful feedback during the drafting process.

\printbibliography
\end{document}